\documentclass[12pt]{article}
\usepackage[cp1251]{inputenc}
\usepackage[english,russian,ukrainian]{babel}
\usepackage{latexsym,amsfonts,amsthm,amssymb,amsmath}
\usepackage{graphicx,graphics,hhline}
\usepackage{euscript}
\usepackage[all]{xy}
\theoremstyle{plain}
\newtheorem{theorem}{Theorem}
\newtheorem{proposal}{Proposition}
\newtheorem{lemm}{Lemma}
\newtheorem{corol}{Corollary}
\theoremstyle{definition}
\newtheorem{definition}{Definition}
\newtheorem{examp}{Example}

\begin{document}

\noindent {\large \bf ON SEMICONVEXITY \\ \\ AND WEAK SEMICONVEXITY }

\large

\vskip 1cm

\noindent UDC 514.172

\noindent MSC 32F17, 52A30

\vskip 1cm

\noindent{\bf T.M.~Osipchuk (Т.М.~Осіпчук)},  osipchuk@imath.kiev.ua
\vskip 2mm

\noindent Institute of Mathematics NAS of Ukraine

\vskip 1cm

{\small
\noindent {\bf Abstract.} Properties of two classes of generally convex sets in the space  $\mathbb{R}^n$, called $m$-semiconvex and weakly $m$-semiconvex, $1 \le m<n$, are investigated in the present work.   In particular, it is established that an open set with smooth boundary in the plan  which is weakly $1$-semiconvex but not $1$-semiconvex  consists minimum of four simply connected components.

\noindent {\bf Анотація.}
В роботі досліджуються  властивості двох класів узагальнено опуклик множин у просторі  $\mathbb{R}^n$, які називаються $m$-напівопуклими та слабко $m$-напівопуклими, $1 \le m<n$. Зокрема, встановлено,
що відкрита множина з гладкою межею на площині, яка є слабко $1$-напівопуклою, але не $1$-напівопуклою, складається мінімум з чотирьох однозв'язних компонент.

\noindent {\bf Аннотация.} { В работе исследуются свойства двух классов обобщенно выпуклых множеств в пространстве
  $\mathbb{R}^n$, которые называются $m$-полувыпуклыми и слабо $m$-полувыпуклыми, $1 \le m<n$.
Вчастности, установлено, что открытое множество с гладкой границей на плоскости, являющееся слабо $1$-полувыпуклым но не $1$-полувыпуклым, состоит минимум из четырех односвязных компонент.
  }

}

\normalsize

\section{Introduction}

The class of $m$ - semiconvex sets is one of the classes of generally convex sets.  The semiconvexity notion was proposed by Yu.~Zeliskii \cite{Zel1} and it was used in the formulation of a shadow problem generalization. The shadow problem was proposed by G. Khudaiberganov  \cite{Hud} in 1982  and is the following: {\it To find the minimal number of open (closed) balls in the space $\mathbb{R}^n$ that are pairwise disjoint, centered on a sphere $S^{n-1}$ (see \textup{\cite{Roz1_1}}),  do not contain the sphere center, and such that any straight line passing through the sphere center intersects at least one of the balls}. To formulate the generalized shadow problem, first, let us give the following definitions which we also use in our investigation.

Any $m$-dimensional affine subspace of Euclidean space  $\mathbb{R}^n$, $m< n$, is called an {\it$m$-dimensional plane}.
\begin{definition}\label{def0}
One of two parts of an $m$ - dimensional plane of the space $\mathbb{R}^n$, $n\ge 2$, into which it is divided by its any $(m-1)$ - dimensional plane, is said to be an {\it $m$ - dimensional half-plane}.
\end{definition}
An $m$ - dimensional  half-plane can be either open or closed.  An {\it open} $m$ - dimensional  half-plane is one of two open sets produced by the subtraction of an $(m-1)$ - dimensional plane from an $m$ - dimensional plane. A {\it  closed} $m$ - dimensional  half-plane is the union of an open $m$ - dimensional half-plane and the $(m-1)$ - dimensional plane that defines it.

 For instance, the $1$ - dimensional half-plane is a ray, the $2$ - dimensional half-plane is a half-plane, etc.

\begin{definition}\label{def1} (\cite{Zel2}) A set $E\subset\mathbb{R}^n$ is called {\it m-semiconvex with respect to a point} $x\in \mathbb{R}^n\setminus E$,  $1\le m<n$, if there exists
an m-dimensional half-plane $L$ such that $x\in L$ and $L\cap E=\emptyset$.
\end{definition}

\begin{definition}\label{def2} (\cite{Zel2}) A set $E\subset\mathbb{R}^n$ is called {\it m-semiconvex },  $1\le m<n$, if it is m-semiconvex with respect to every point $x\in \mathbb{R}^n\setminus E$.
\end{definition}

One can easily see that both definitions satisfy the axiom of convexity: The intersection of
each subfamily of these sets also satisfies the definition. Thus, for any set $E\subset\mathbb{R}^n$ we can consider the minimal $m$-semiconvex set containing $E$. This set is called the {\it $m$-semiconvex hull of set  $E$}.

The generalized shadow problem is {\it To find the minimum number of pairwise disjoint closed (open) balls in $\mathbb{R}^n$ (centered on the sphere $S^{n-1}$ and whose radii are smaller than the radius of the sphere)  such that any
ray starting at the center of the sphere necessarily intersects at least one of these balls}.

In the terms of $m$-semiconvexity this problem can be reformulated as follows: {\it What is the minimum number of pairwise disjoint closed (open) balls  in $\mathbb{R}^n$ whose  centers are located on a sphere $S^{n-1}$ and the radii are smaller than the radius of this sphere such that the center of the sphere belongs to the $1$-semiconvex hull of the family of balls?}

In the paper \cite{Zel1} the problem is solved as $n=2$. And only the sufficient number of the balls is indicated as $n=3$.

In the 60's L.~Aizenberg and A.~Martineau proposed their notions of a linearly convex set in the multi-dimensional  complex space $\mathbb{C}^n$. The first author considered domains and their closures and used boundary points of the domains in his definition \cite{Aiz3}-\cite{Aiz1}. The second one used all points of the addition to a set of the space $\mathbb{C}^n$ \cite{Martino2}.  If one uses these definitions not only for domains and compact sets, then Aizenberg's definition isolates one connected component of a set linearly convex in the sense of Martineau.

Guided by similar reflections,  Yu.~Zeliskii suggested to distinguish $m$-semiconvex and weakly $m$-semiconvex sets and obtained the following results.

\begin{definition}\label{def3} (\cite{Zel3}) An open set  $G\subset\mathbb{R}^n$ is called {\it
weakly $m$-semiconvex},  $1\le m<n$, if it is $m$-semiconvex for any point  $x\in \partial G$. A set $E\subset\mathbb{R}^n$ is called {\it weakly $m$-semiconvex} if it can be approximated from the outside by a family of open weakly $m$-semiconvex sets.
\end{definition}

\begin{theorem}\label{lemm1} \textup{(\cite{Zel3} )} Let a set $E\subset \mathbb{R}^2$ be weakly $1$-semiconvex and not $1$-semiconvex. Then set $E$ is disconnected.
\end{theorem}

In \cite{Zel3} there was also made the assumption that a weakly $1$-semiconvex and not $1$-semiconvex set consists of not less than three components. This proposition was proved in \cite{Dak9}.

\begin{theorem}\label{lemm11} \textup{(\cite{Dak9} )} Let a set $E\subset \mathbb{R}^2$ be weakly $1$-semiconvex and not $1$-semiconvex. Then $E$ consists of not less than three components.
\end{theorem}

The present work proceeds  the research of Yu. Zelinskii  and contains the investigation of properties of $1$-semiconvex and weakly $1$-semiconvex sets with smooth boundary of the plane.

Except  open and closed sets,  a set that contains a subset $Q$ of the set $B$ of its boundary points  and does not contain $B\setminus Q$ is also used in this work.  The expression "neither open nor closed set"\/ will stand for such a set.

\section{Main results}

We consider a class of weakly $m$-semiconvex sets by the following definition:

\begin{definition}\label{def4} A set  $E\subset\mathbb{R}^n$ is called {\it
weakly $m$-semiconvex},  $1\le m<n$, if for any point  $x\in \partial E$ there exists  an open $m$-dimensional
half-plane $L$ such that  $x\in \overline{L}$  and  $L\cap E=\emptyset$.
\end{definition}

Classes of open sets, weakly $m$-semiconvex by Definition \ref{def3} and Definition \ref{def4}, coincide. However, a set containing  closed or  neither open nor closed connected components which is not weakly $m$-semiconvex by Definition \ref{def3} can belong to the class of weakly $m$-semiconvex sets with respect to Definition \ref{def3}, Figure 4, 6.  And conversely, there are sets weakly $m$-semiconvex by Definition \ref{def4} but not weakly $m$-semiconvex with respect to Definition \ref{def4}, as illustrated in Figure 1.

Everywhere further we consider weakly $m$-semiconvex sets with respect to Definition \ref{def4}. We say that {\it there exists (can be drawn) an $m$-dimensional half-plane starting at a point and satisfying some conditions} if there exists an open $m$-dimensional half-plane satisfying the conditions and which boundary contains the point.

Let us provide a number of accessory propositions.

\begin{proposal}\label{pro1}  If an open set $E\subset\mathbb{R}^n$ is $m$-semiconvex, then it  is weakly $m$-semiconvex.
\end{proposal}

It fails for closed and neither open nor closed sets. Figure 1 shows examples of sets that are $1$-semiconvex but not weakly $1$-semiconvex.

\begin{center}
\includegraphics[width=10 cm]{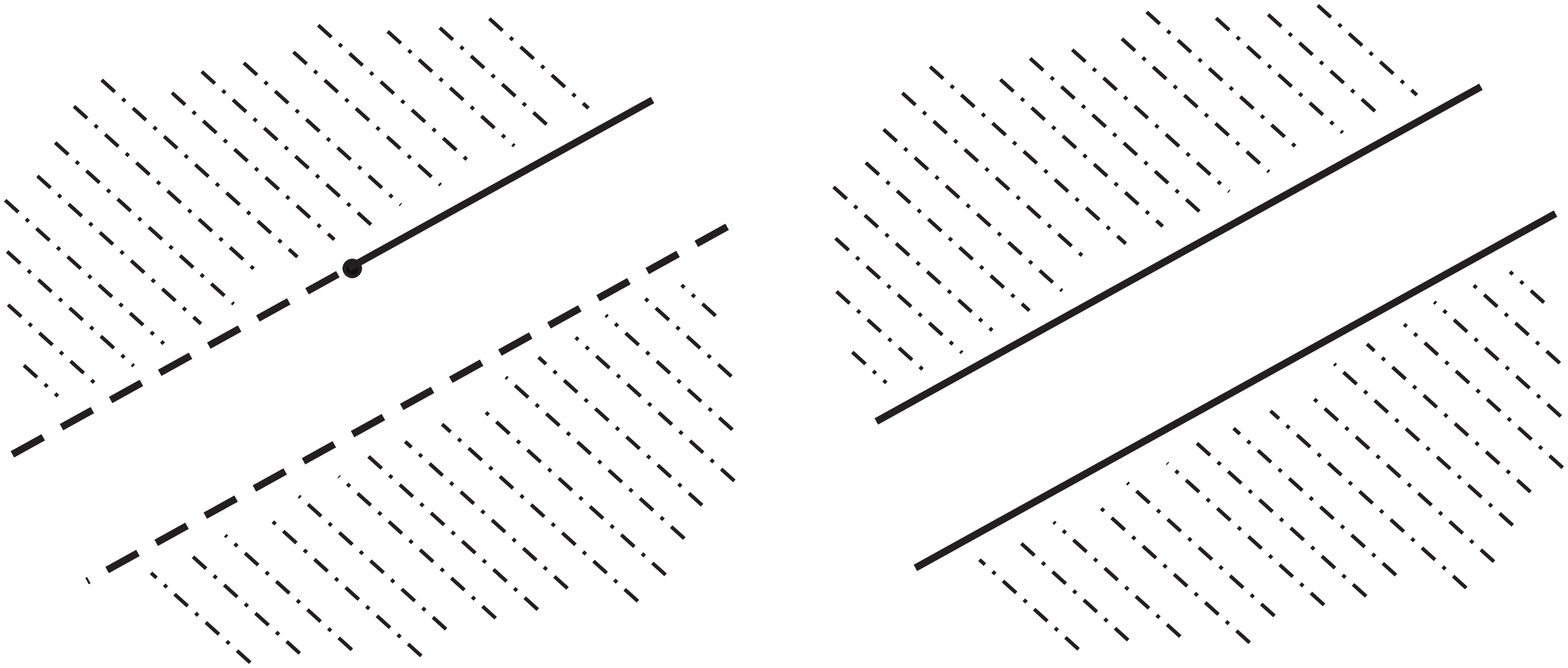}
\end{center}
\small\begin{center} Fig. 1 \end{center}
\normalsize

\begin{proposal}\label{pro2}  A weakly $(n-1)$-semiconvex \emph{ (}$(n-1)$-semiconvex \emph{)} set $E\subset\mathbb{R}^n$ consists of simply connected components.
\end{proposal}
\begin{proposal}\label{pro3}  There exist weakly $1$-semiconvex sets in the plane which are not $1$-semiconvex.
\end{proposal}

\begin{examp}\label{examp1}
Figure 2 a) shows the set $E$ that consists of four open rectangles with some common tangents. For any boundary point of set $E$ there exists a ray that does not intersect  $E$ but for any point of the interior of rhombus $ABCD$ such a ray can not be found.  In Figure 2 b), the set consisting of three open components is  weakly $1$-semiconvex, but  for any point of the interior of triangle $ABC$ we can not find a ray that does not intersect the set.
\end{examp}
 \begin{center}
\includegraphics[width=12 cm]{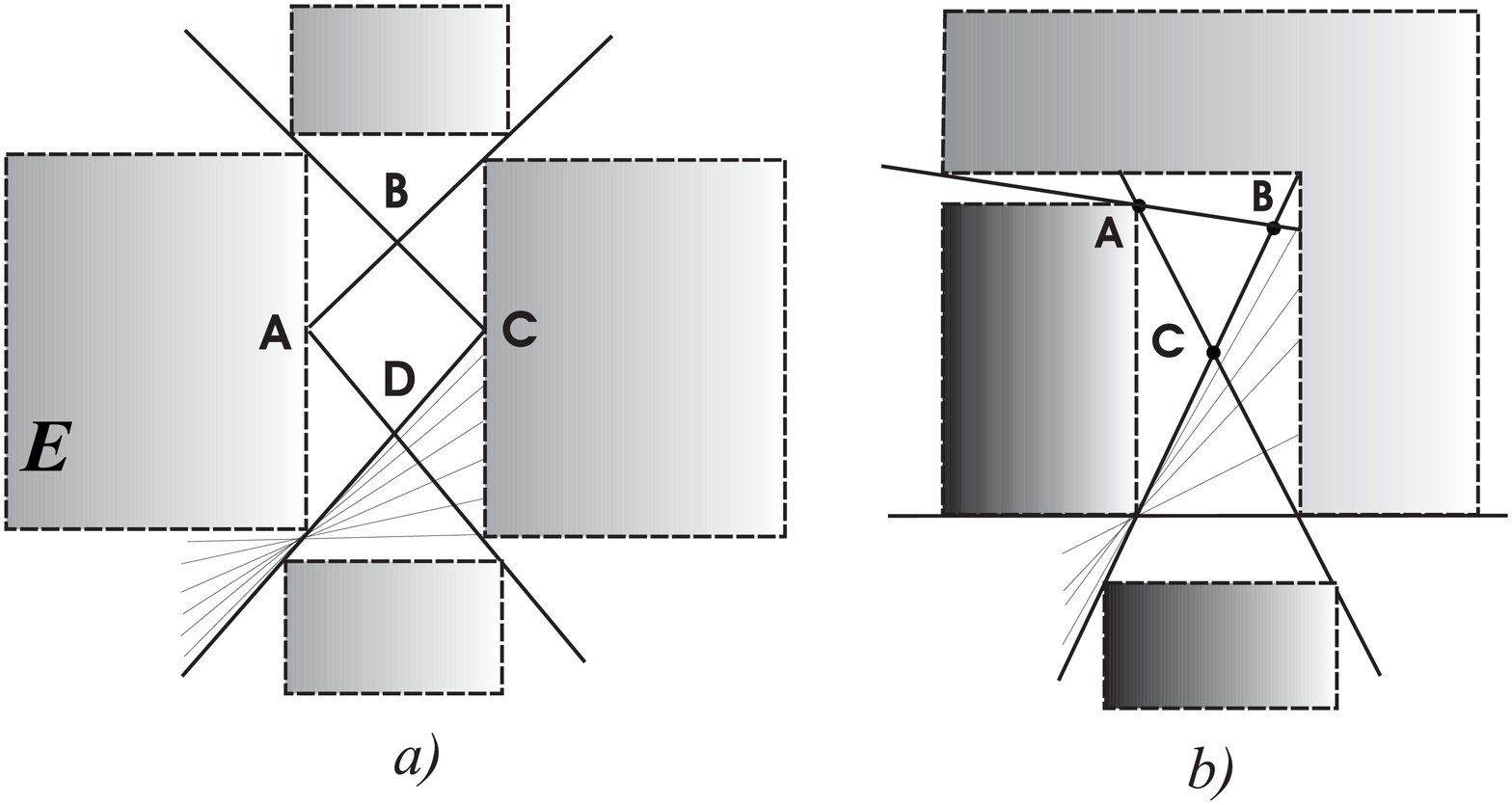}
\end{center}
\small\begin{center} Fig. 2 \end{center}
\normalsize

\begin{theorem}\label{lemm1} Let a set $E\subset \mathbb{R}^2$ be weakly $1$-semiconvex and not $1$-semiconvex. Then set $E$ is disconnected.
\end{theorem}

\noindent{\bf Proof.}  Since set $E$ is weakly $1$-semiconvex and not $1$-semiconvex,  there is a point $x\in \mathbb{R}^2\setminus \overline{E}$ such that any ray starting at point $x$ intersects set $E$. We now fix two complementary (lying on the same straight line) rays $s_1$, $s_2$ starting at point $x$. They intersect $\partial E$ at some points $x', x''$ that are nearest to $x$ along rays $s_1$, $s_2$ respectively. But,  since $E$ is a weakly $1$-semiconvex set, it follows that for points $x', x''$\/ there exist rays that do not cross $E$.  The polygonal chain containing rays $s_1$, $s_2$, with vertexes at points $x',  x''$,  cuts then the plane into two parts if it is simple (without self-intersections) (case 1) or into three parts if the polygonal chain is self-intersecting (case 2).

If in case 1 the whole set $E$ is contained in one part of the plane,  then for point $x$ we can always draw a ray  lying in the other part. Such a ray does not intersect $E$, which contradicts the lemma conditions.

In case 2 set $E$ can not be contained in one part of the plane too. Indeed, if the whole set $E$ is contained in part I (see Figure 3), then ray $s$ starting at the point $x$ and passing through the point of intersection of  rays $s_1$, $s_2$ is a ray that does not intersect $E$. If set $E$ is in one of the other parts, then the ray complementary to  $s$  can be the required ray.

 \begin{center}
\includegraphics[width=7 cm]{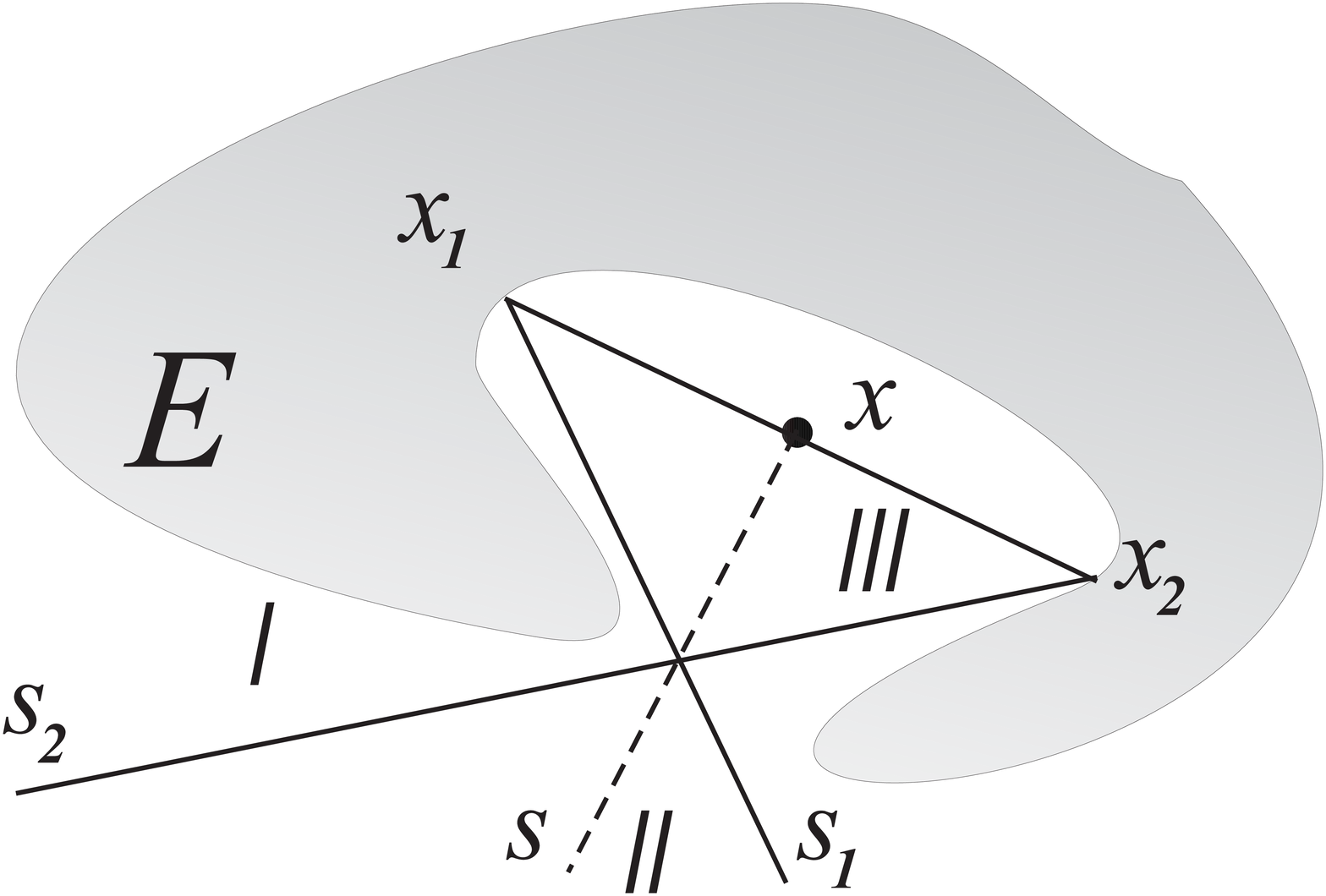}
\end{center}
\small\begin{center} Fig. 3 \end{center}
\normalsize

The idea of the proof of Lemma 1 belongs to Yu. Zelinskii (\cite{Zel3}).

Theorem \ref{lemm11} is true in the sense of Definition \ref{def3}. But  there are sets that are weakly $1$-semiconvex but not $1$-semiconvex by Definition \ref{def4} and consisting of two components,  as the following example shows.

\begin{examp}
Let a set $E\subset\mathbb{R}^2$ consist of two closed components with common tangent line and suppose $E$ is not smooth at the points $A$, $B$, Figure 4. It is weakly $1$-semiconvex but not $1$-semiconvex at all points of the open segment $(A,B)$.
\end{examp}

\begin{center}
\includegraphics[width=5 cm]{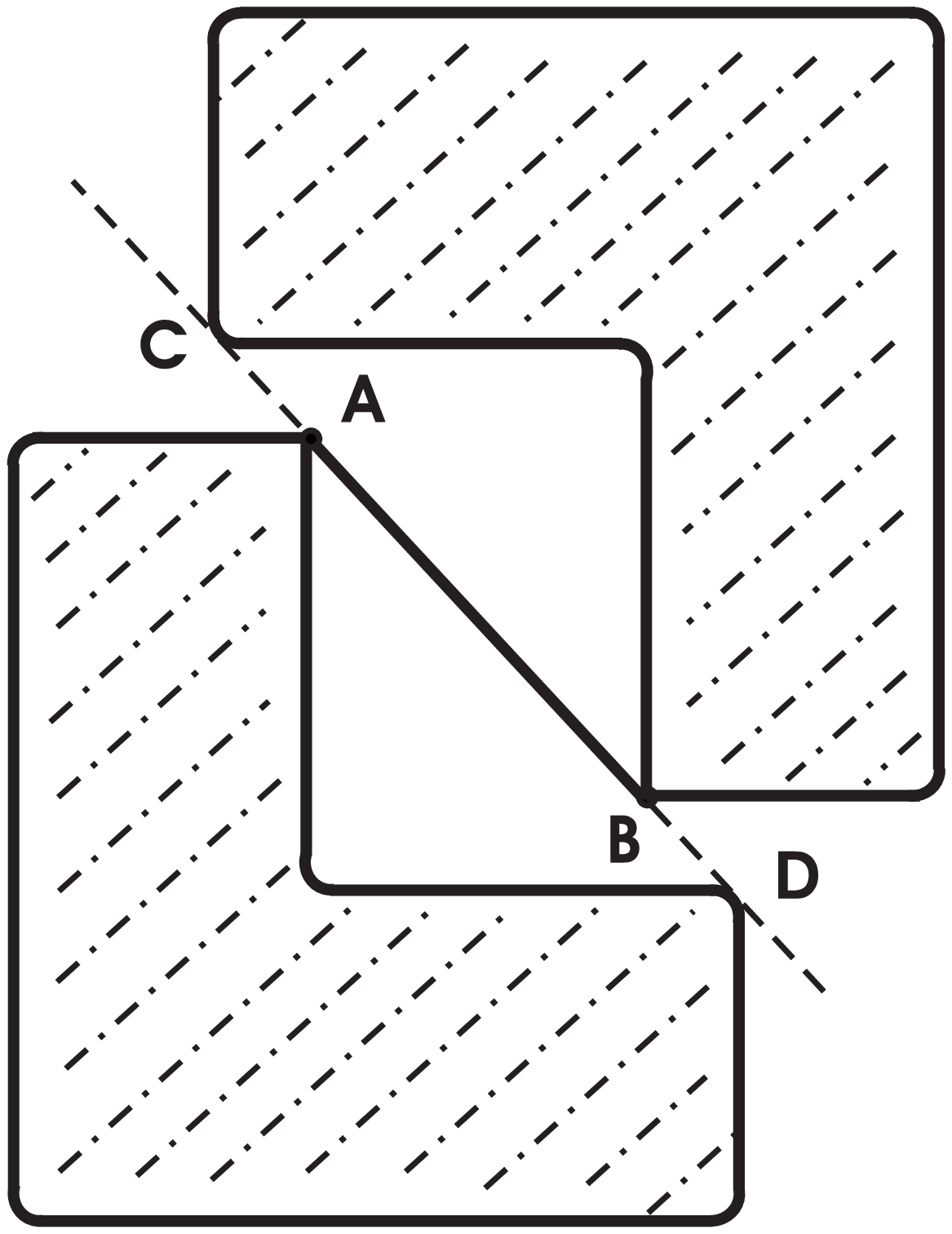}
\end{center}
\small\begin{center} Fig. 4 \end{center}
\normalsize

Here we propose our own result, similar to  Theorem \ref{lemm11}, but for sets with smooth boundary.

\begin{theorem}\label{theor2}  Let  $E\subset \mathbb{R}^2$ be an open (closed),  weakly $1$-semiconvex and not $1$-semiconvex set with smooth boundary. Then set $E$ consists of more than two simply connected components.
\end{theorem}

\noindent{\bf Proof.}  Set $E$ is disconnected by Theorem \ref{lemm1} and all its components are simply connected by Proposition \ref{pro2}. Suppose set $E$ consists of two simply connected components $E_1$, $E_2$. First, we consider the case when $E$  is  open.

Since set $E$ is not $1$-semiconvex and weakly $1$-semiconvex, it follows that  there is a point $x\in \mathbb{R}^2\setminus \overline{E}$ such that any ray starting at $x$ intersects set $E$ and its boundary $\partial E$ as well. Suppose a ray $l_1$ starting at point $x$ and intersecting set $E$ crosses $\partial E$ at some point, nearest to $x$ along the ray. Since set $E$ is  weakly $1$-semiconvex, there is a ray $l$ starting at this point and not intersecting $E$.

Among all rays starting at point $x$ and crossing ray $l$, there exist such rays that intersect $E$ before they intersect $l$. A case  when all these rays intersect $E$ after they intersect $l$ is possible only when $E$ is unbounded. But in this case there exists the ray that is parallel to the ray $l$ and does not intersect $E$, which contradicts the conditions of the theorem.

For definiteness, suppose the rays from the previous paragraph initially  intersect  component  $E_1$. Then from point $x$ we can draw a ray $s_1$ that is tangent to $E_1$ at some point $y_1$ and does not intersect $E_1$. By the way, the extreme positions of  ray $s_1$ is  ray $l_1$. Suppose  $s_1$ intersects $l$ at some point  $z_1$.

For point $y_1$ we can draw a ray that does not intersect $E$,  by the conditions. So we can draw the ray that is tangent to $E_1$ and does not intersect $E$.  With that, such a ray is to pass above the straight line $s$ that contains  ray $s_1$, because  $s$ is tangent to the boundary of $E_1$ and intersects set $E$ (see Figure 5). Here and everywhere below the expression "to pass above the straight line $s$"\/ means to pass in the half-plane with respect to $s$  which does not contain $E_1$ in some neighborhood of point $y_1$. Such neighborhood exists according to the constructions in the previous paragraph.

\begin{center}
\includegraphics[width=10 cm]{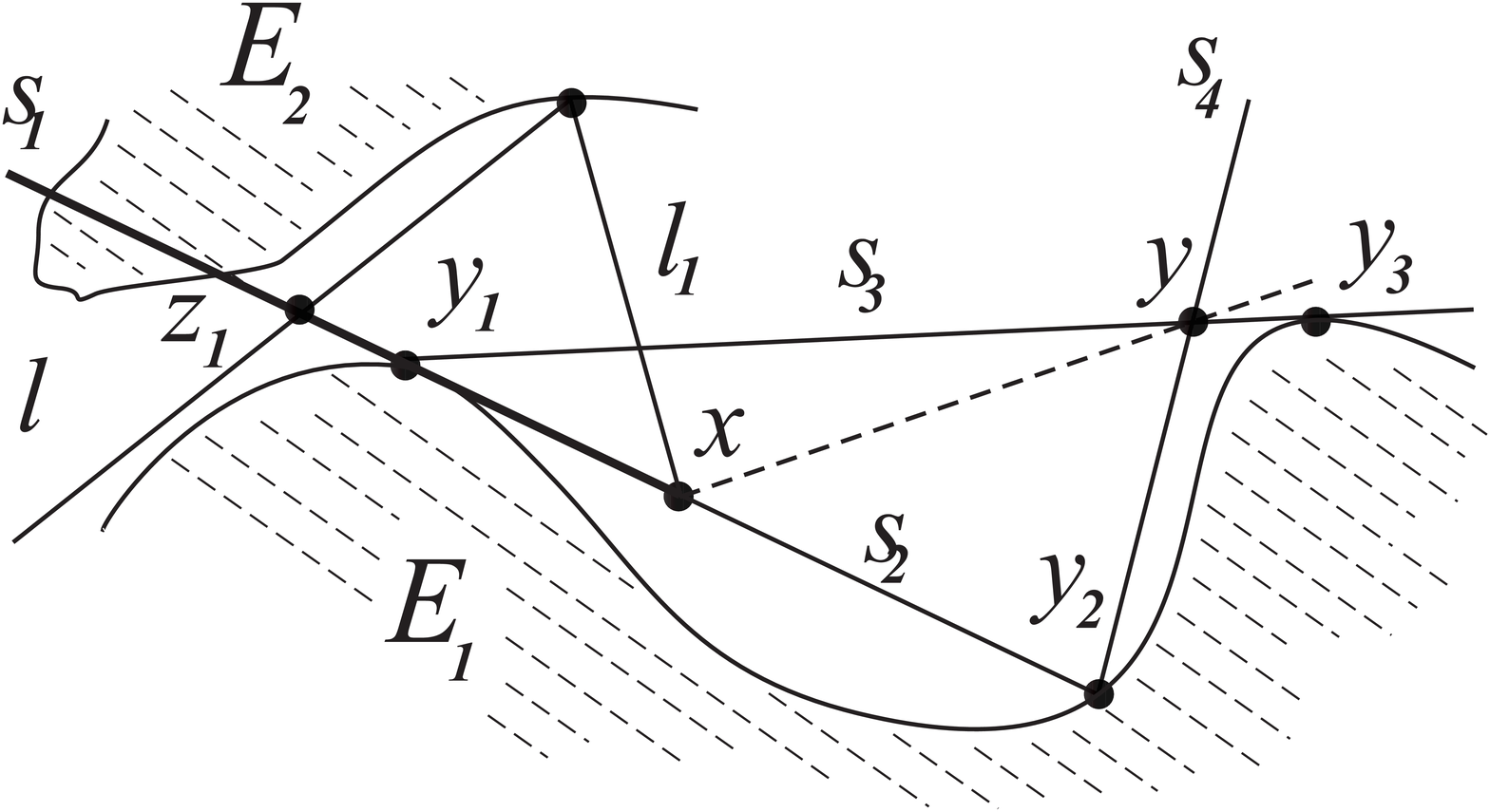}
\end{center}
\small\begin{center} Fig. 5 \end{center}
\normalsize

Let a ray $s_3$, starting at the point $y_1$, be tangent to the component $E_1$ at some point $y_3$ that is above the straight line $s$, according to the previous explanations. Points  $y_1$, $y_3$ belong to the boundary of the same component $ E_1$. Ray $s_1$ does not intersect the part of $\partial E_1$ connecting the points, but the ray $s_2$, complement to  $s_1$, intersects $\partial E_1$ at some point $y_2$. For $y_2$ we can draw  a ray $s_4$ that does not intersect set  $E$ but it should  intersect the segment $y_1y_3$ at some point $y$.

According to the previous constructions, the polygonal chain $lz_1y_1s_3$ is simple (without self-intersections) and it cuts the plane into two parts, each of which contains a component of set $E$. With that, the domain that is bounded by the part of $\partial E_1$, connecting points $y_1$, $y_3$, on one side and by the segment $y_1y_3$ on the other one, does not contain points of $E$. Thus,  we can draw the ray, starting at point $x$ and passing through the point $y$, which does not intersect  $E$. But this is a contradiction with the theorem conditions.

In the case of closed set $E$ all statements remain the same as in the case of open set $E$ with the following small addition.  The ray $s_1$ already intersects  $E$ at the point $y_1$. But still it should intersect $E$ after crossing $l$.  In the opposite case, for point $x$ we can draw a ray that lies between  $s_1$ and the ray tangent to  $E$ behind $l$ and does not intersect $E$, which contradicts the theorem conditions.

Thus, $E$ consists of more than two simply connected components. Theorem \ref{theor2} is proved.

\begin{examp}
The example of weakly $1$-semiconvex and not $1$-semiconvex set with smooth boundary which consists of two simply connected components is given in Figure 6. In  case a), the set $E$, among all its boundary points, contains only points $A$ and $B$. In case b) $E$ contains all its boundary points except points $C$ and $D$. In both cases, for any point of the interval $AB$ it is not possible to draw a ray that does not intersect the sets. In general, if one cuts a subset of the set of all boundary points of $E$ which contains points $C$, $D$ and does not contain points  $A$,  $B$ from the closure of  $E$, then the obtained set is also weakly $1$-semiconvex and not $1$-semiconvex.
\end{examp}

\begin{center}
\includegraphics[width=10 cm]{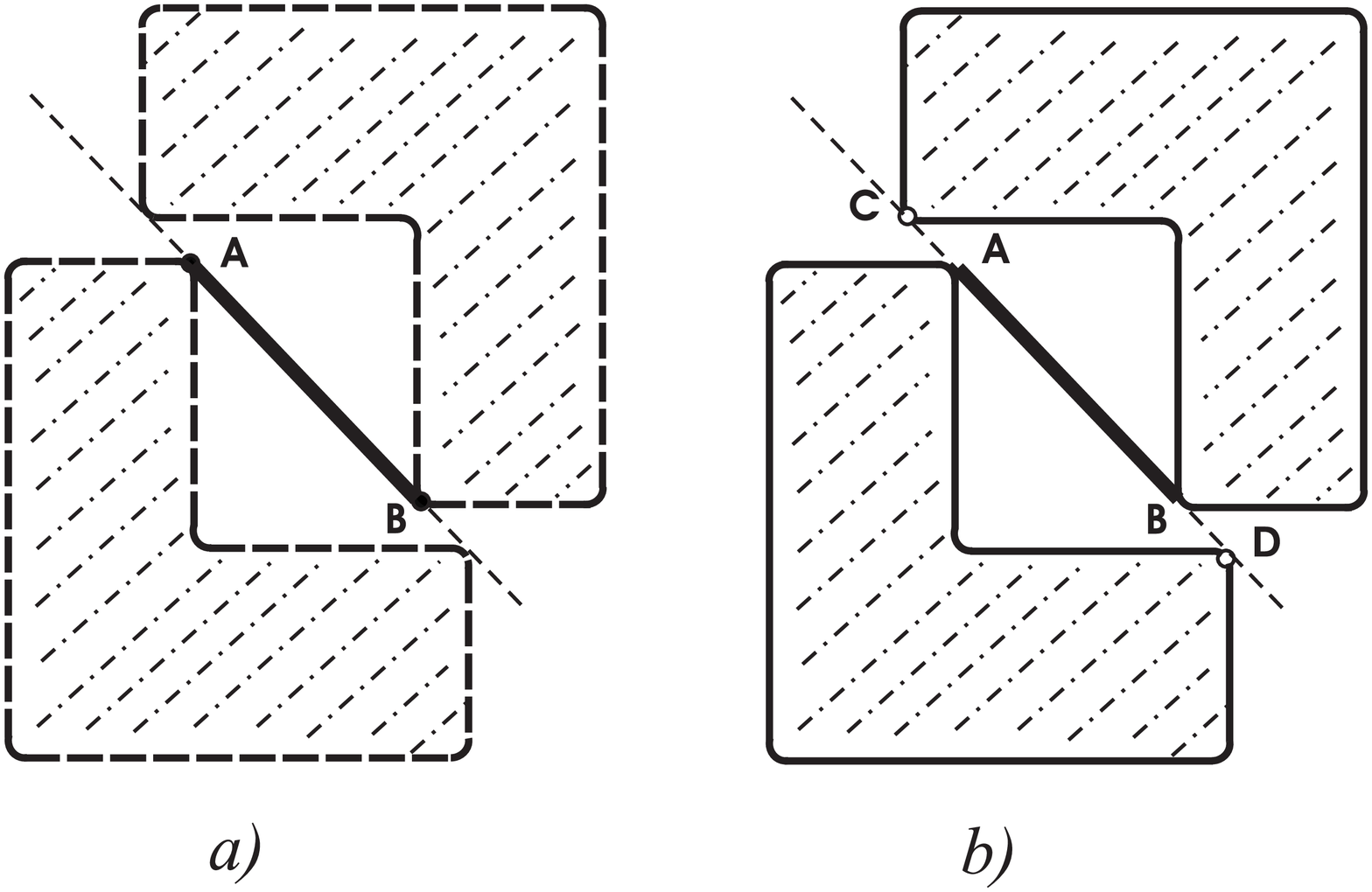}
\end{center}
\small\begin{center} Fig. 6 \end{center}
\normalsize

\begin{examp}\label{examp3}
The system of four open (closed) balls with common tangent lines, illustrated in Figure 7a (7b), is an example of open (closed) set with smooth boundary which is weakly $1$-semiconvex, not $1$-semiconvex, and consists of four simply connected components. Herewith, it is easy to see that in the case of open balls for any point of the interior of the rhombus $ABCD$ there is no ray that intersects the set (Figure 7a). And in the case of closed balls this holds for the points of the closure of  rhombus $ABCD$ (Figure 7b).
\end{examp}
\begin{center}
\includegraphics[width=15 cm]{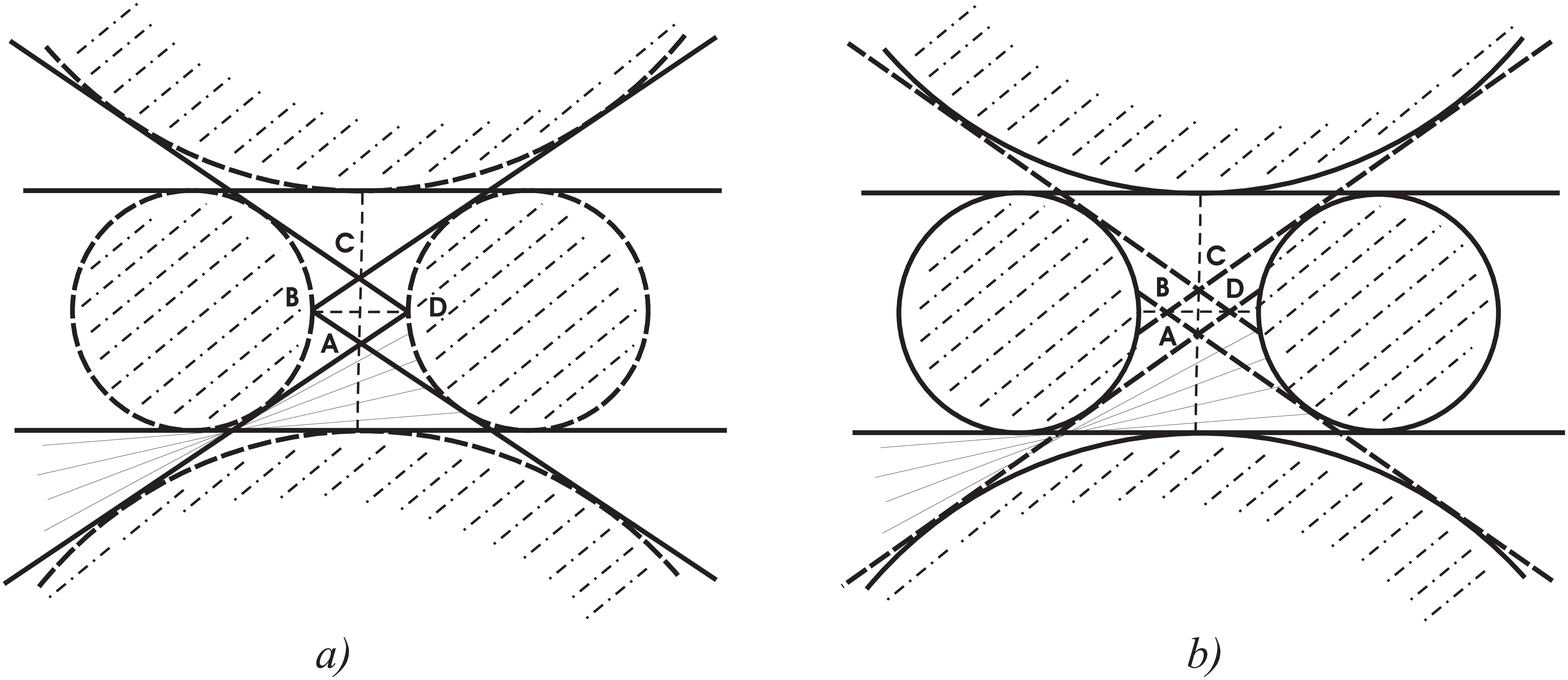}
\end{center}
\small\begin{center} Fig. 7 \end{center}

\normalsize

The following result is the direct corollary from Theorem \ref{lemm1}.

\begin{lemm}\label{corol1}
Let a set $E\subset \mathbb{R}^2$ be simply connected and weekly $1$-semiconvex. Then $E$ is $1$-semiconvex.
\end{lemm}

\begin{definition}\label{def5} A point  $x\in \mathbb{R}^n\setminus A$ is called {\it
point of  $m$-nonsemiconvexity of a set $A\subset\mathbb{R}^n$} if there is no open $m$-dimensional half-plan starting at $x$ and not intersecting  $A$.
\end{definition}

\begin{definition}\label{def6} A ray  $l$ starting at a point $x\in \mathbb{R}^n\setminus\overline{A}$ is  {\it
supporting} for a set $A\subset \mathbb{R}^n$ at a point $a\in \partial A$ if $a\in l\cap \partial A\ne\emptyset$  and $l\cap \mathrm{Int}\, A=\emptyset$.
\end{definition}

\begin{lemm}\label{lemm2} Let an open (closed) set $E\subset \mathbb{R}^2$ be simply connected, weakly $1$-semiconvex, and have smooth boundary. Then there exist two and only two supporting rays of $E$ starting at a point $x\in \mathbb{R}^2\setminus \overline{E}$.
\end{lemm}

\noindent{\bf Proof.}  The set $E$ is $1$-semiconvex by Corollary \ref{corol1}. So, for any point $x\in \mathbb{R}^2\setminus \overline{E}$ we can draw a ray that does not intersect $E$. By continuity  and since $E$ has smooth boundary, we can draw at least one tangent to $\partial E$ ray starting at $x$ that does not intersect $\mathrm{Int}\, E$.

Suppose we have a unique supporting ray for a fixed point $x$. This is not possible for closed sets. Let set $E$ be open. Without loss of generality we suppose that there are two possible cases: $\partial E$ and supporting ray  have one common point $x_1$;  $\partial E$ and supporting ray  have two common points $x_1$, $x_2$, with that, $x_1$ is closer to $x$ than $x_2$. The complimentary ray to the supporting one intersects $\partial E$ at some point $x_0$. Since $\partial E$ is smooth, points of both open parts of $\partial E$ between  $x_1$, $x_0$ are $1$-nonsemiconvexity points in the first case (Figure 8a)), while  points of open part of $\partial E$ between  $x_1$, $x_0$ are $1$-nonsemiconvexity points in the second one (Figure 8b)). In both cases this is a contradiction with  $1$-semiconvexity of $E$.
\begin{center}
\includegraphics[width=10 cm]{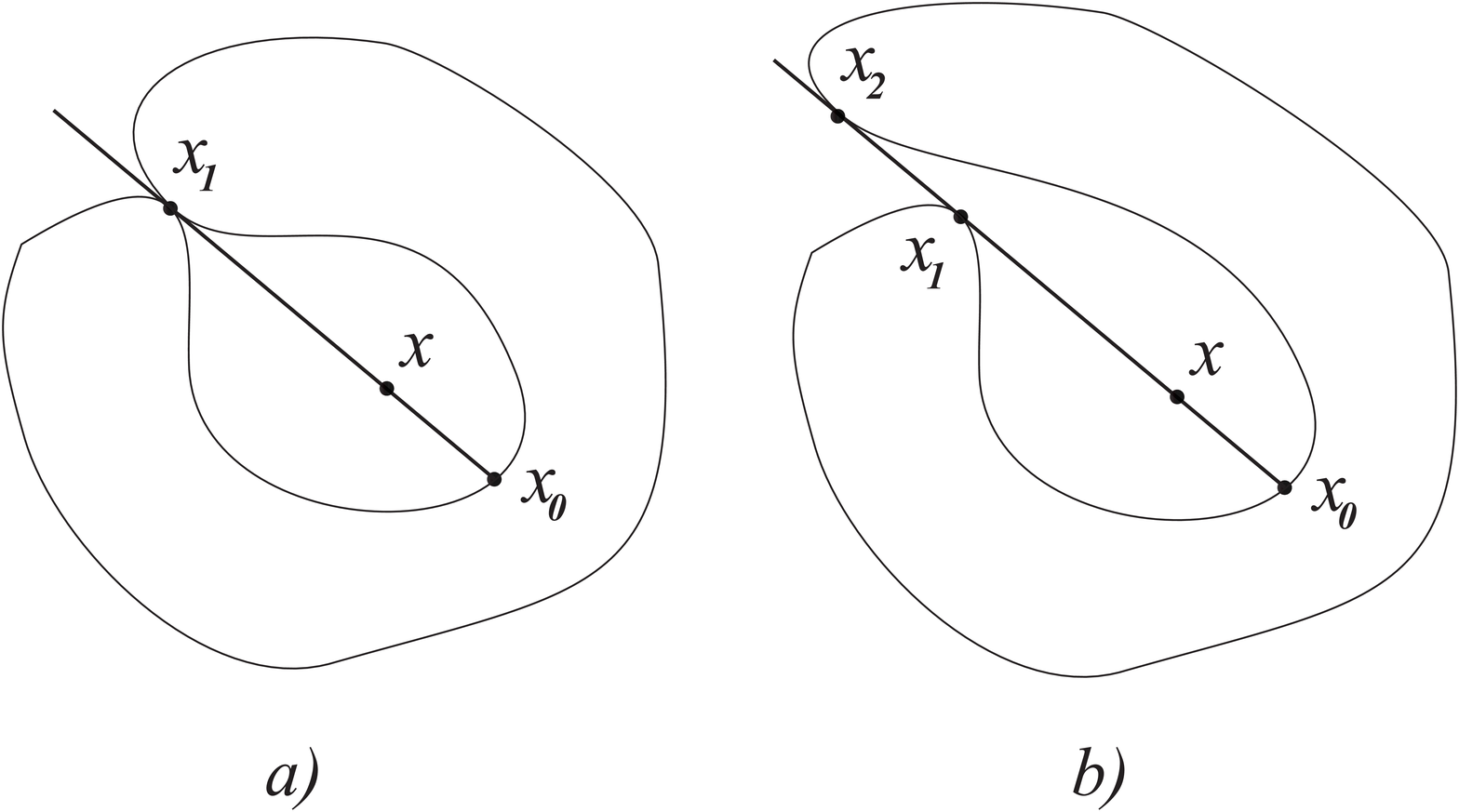}
\end{center}
\small\begin{center} Рис. 8 \end{center}
\normalsize

Let us  assume that there exist $n$, $n\ge 3$, different supporting rays starting at the point $x\in \mathbb{R}^2\setminus \overline{E}$. So they cut the plane into $n$ parts and, by the construction,  $n-1$ of them should contain $E$, which contradicts the fact that  $E$  is connected.

The lemma is proved.

\vskip 2mm

The open set on Figure 8 a) with boundary that is nonsmooth  at the point $x_1$ can have a unique supporting ray starting at the point $x$.

\vskip 2mm

We say that {\it a set $A\subset \mathbb{R}^n$ is projected from a point $x\in\mathbb{R}^n$ on a set $B\subset \mathbb{R}^n$} if any ray, starting at point $x$ and intersecting $A$,  intersects $B$ as well.

\vskip 2mm

\begin{lemm}\label{lemm3} Let an open set $E\subset \mathbb{R}^2$ be weakly $1$-semiconvex but not $1$-semiconvex and consist of three components. Then none of its components is projected on the others from a point of $1$-nonsemiconvexity of $E$.
\end{lemm}

\noindent{\bf Proof.}  Let $E$ consist of three components $E_1$, $E_2$, $E_3$ and  $x\in \mathbb{R}^2\setminus \overline{E}$ be a point of $1$-nonsemiconvexity of $E$.  Without loss of generality, suppose $E_1$ is projected from $x$  on at least one of the other components, Figure 9. Then the set, consisting only of components $E_2$, $E_3$, is weakly $1$-semiconvex and not $1$-semiconvex, which contradicts Theorem \ref{lemm11}.
\begin{center}
\includegraphics[width=12 cm]{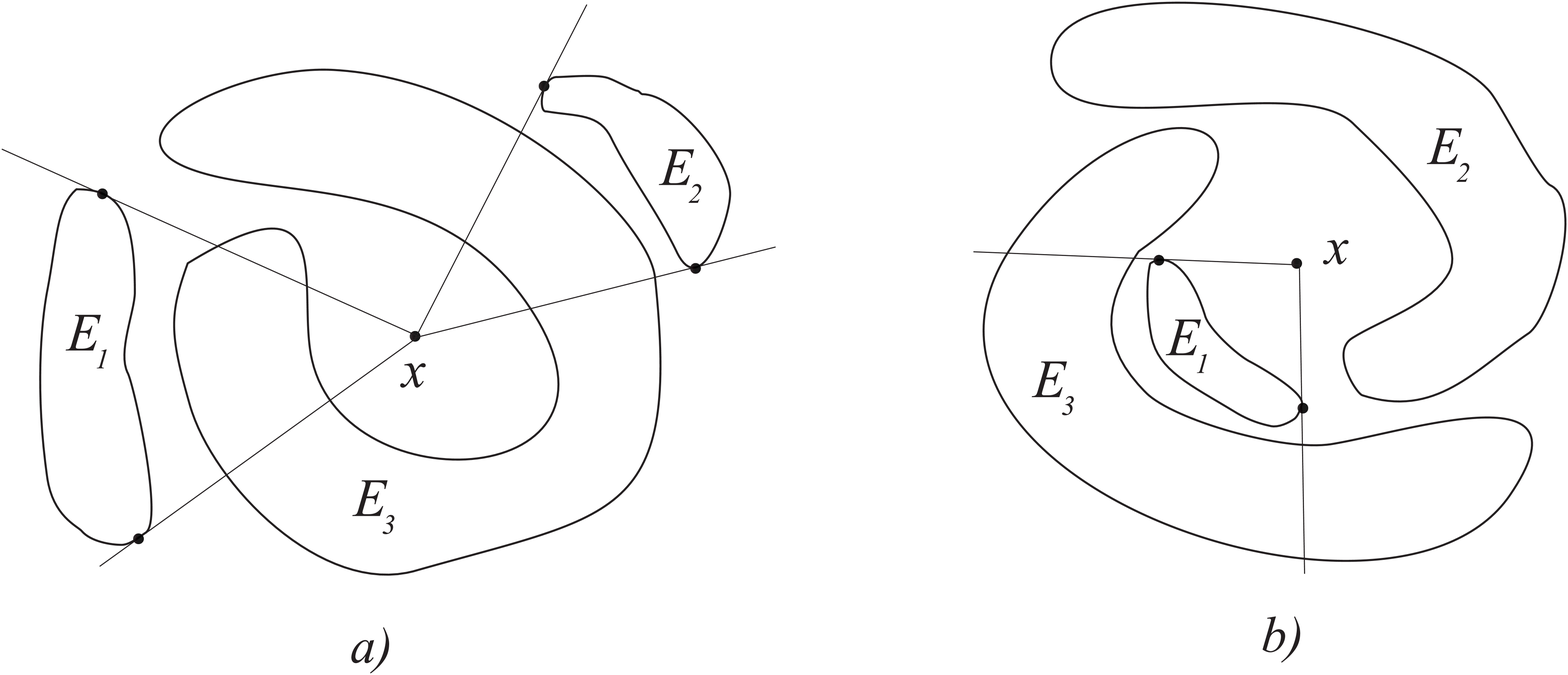}
\end{center}
\small\begin{center} Fig. 9 \end{center}
\normalsize

The lemma is proved.
\vskip 2mm

Let sets $A, B\subset \mathbb{R}^n$ be given. Let $l$ be a supporting ray of $A$ starting at a point $x\in\mathbb{R}^n$ and point $a\in l\cap \partial A$. The ray $l$ is called {\it inner supporting ray of the set $A\cup B$} if $l\cap B\ne\emptyset$ and the distance between points $x$ and $a$  is less than a distance between point $x$ and any point $b\in l\cap B$.

\vskip 2mm

\begin{lemm}\label{lemm4} Let a set $E\subset \mathbb{R}^2$ be weakly $1$-semiconvex but not $1$-semiconvex. Then $E$ has at least one inner supporting ray.
\end{lemm}

\noindent{\bf Proof.} Set $E$ consists of $n$, $2<n<\infty$, simply connected components $E_j$, $j=\{\overline{1,n}\}$, by Theorem \ref{lemm11}.  Since $E$ is not $1$-semiconvex, there exists a $1$-nonsemiconvexity point $x\in \mathbb{R}^2\setminus \overline{E}$. Let us draw a ray $\xi$ starting at point $x$. It intersects $\partial E$ at some point $y$ that is nearest to $x$ along the ray. Since $E$ is weakly $1$-semiconvex, there exists a ray $\gamma$ starting at $y$ and not intersecting $E$. Moreover, $\gamma$ does not lie on the straight line that contains $\xi$.

Let us rotate ray $\xi$ round $x$ in the half-plane where $\gamma$ lies.  In other words, we choose the polar coordinate system where point $x$ is the pole,  ray $\xi$ is the polar axis and a positive angular coordinate is determined by the angle $\varphi$ between the ray $\xi=\xi(0)$ and a ray $\xi(\varphi)$, $0<\varphi<\pi$, that intersects ray $\gamma$.

So, among all rays starting at $x$ and intersecting $\gamma$, there are those that intersect component $E_1$, for definiteness,  before they intersect $\gamma$. A case  when all these rays intersect $E$ after they intersect $\gamma$ is possible only when $E$ is unbounded. But in this case there exists the ray that is parallel to ray $\gamma$ and does not intersect $E$, which contradicts the lemma conditions. Thus, by continuity, there exists the ray $\xi(\varphi_0)$, $0\le\varphi_0<\pi$, such that $\xi(\varphi_0)\cap \partial E_1\ne\emptyset$, $\xi(\varphi_0)\cap \mathrm{Int} E_1=\emptyset$, and for $\varepsilon>0$ small enough rays $\xi(\varphi)$, $\varphi_0<\varphi<\varphi_0+\varepsilon$, intersect $E_1$ before they intersect $\gamma_1$. Ray $\xi(\varphi_0)$ does not intersect component $E_1$ after intersecting $\gamma$, otherwise, the point of  $\partial E_1$ which is nearest to point $x$ along a ray  $\xi(\varphi)$, $\varphi_0<\varphi\le\varphi_0+\pi$,   is a  point of 1-nonsemiconvexity of $E$. Thus, $\xi(\varphi_0)\cap E_2\ne\emptyset$, for definiteness, and $\xi(\varphi_0)$ is an inner supporting ray of $E$.

\vskip 2mm

Further in this work we will often use the constructions of the proof of Lemma \ref{lemm4}.

The main result of the paper is the following

\begin{theorem}\label{theor3}  Let  $E\subset \mathbb{R}^2$ be an open,  weakly $1$-semiconvex, and not $1$-semiconvex set with smooth boundary. Then $E$ consists minimum of four simply connected components.
\end{theorem}

\noindent{\bf Proof.}  The set does not consist of one or two simply connected components by  Theorem \ref{theor2}.  Suppose  $E$ consists of three simply connected components $E_i$,  $i=1,2,3$.

Since $E$ is not $1$-semiconvex, it follows that there exists a $1$-nonsemiconvexity point $x\in \mathbb{R}^2\setminus \overline{E}$ of $E$. Since $E$ is a weakly $1$-semiconvex set, we do not consider points of $\partial E$ as points of $1$-nonsemiconvexity of $E$.

By Lemma \ref{lemm3}, non of the components of $E$ is projected on the others. Then there are three rays $\tau_1$, $\tau_2$, $\tau_3$ starting at $x$ and intersecting a unique component  $E_1$, $E_2$, $E_3$ respectively. That is to say, $\tau_1\cap E\equiv\tau_1\cap E_1$, $\tau_2\cap E\equiv\tau_2\cap E_2$, $\tau_3\cap E\equiv\tau_3\cap E_3$. So, rays $\tau_1$, $\tau_2$, $\tau_3$ cut the plane by three nonempty parts $G_1$, $G_2$, $G_3$.

Without loss of generality, let us consider the closure of domain $G_1$ between the rays $\tau_1$, $\tau_2$.  Since $\tau_3$ intersects component $E_3$ in $\mathbb{R}^2\setminus\overline{G_1}$,  set $\overline{G_1}$ does not contain points of $E_3$. Then $\overline{G_1}$ consists of rays that intersect only $E_1$, only $E_2$, and of those that intersect both $E_1$, $E_2$.  With that, the set of rays, intersecting both $E_1$, $E_2$, is open in $\overline{G_1}$ and its boundary consists of one supporting ray of $E_1$ and one of $E_2$. If one assumes that there is another supporting ray, then it should be supporting for one of the components $E_1$, $E_2$. Without loss of generality, we assume that it is supporting for $E_1$. Since each component of $E$ has exactly two supporting rays, by Lemma \ref{lemm2}, then $E_1$ should be completely contained in  $\overline{G_1}$. But this contradicts the fact that $E_1$ is open and $\tau_1\cap E_1\ne\emptyset$.

  Thus, each closed set $\overline{G_1}$, $\overline{G_2}$, $\overline{G_3}$ contains one and only one inner supporting ray $\xi_i$, $i=1,2,3$ of $E$. Since $E$ has the smooth boundary, each ray $\xi_i$ is also tangent to the boundary of corresponding component at some point $y_i$, $i=1,2,3$.  The rays complementary to $\xi_i$ also intersect $E$ at some points that we denote as $z_i$ respectively.

Since sets $G_1$, $G_2$, $G_3$ have common boundary rays $\tau_j$, $j=1,2,3$, and each inner supporting ray $\xi_i$, $i=1,2,3$, belongs to the respective $\overline{G_i}$, without loss of generality, let us consider the case when two supporting rays coincide with the boundary ray $\tau_3$. Since $\tau_3$ intersects component $E_3$, this ray is supporting for components $E_1$, $E_2$ and is also tangent to their boundaries at points $y_1$, $y_2$ respectively. Points $y_1$, $y_2$ do not coincide; otherwise one can not draw a ray that starts at such a point and does not intersect $E$, which contradicts weakly $1$-semiconvexity of $E$.

Let  point $y_1$ be nearer to point $x$ than to point $y_2$, for definiteness.  Then, let us draw a ray $\gamma_1$ starting at $y_1$ and not intersecting $E$. It lies above the straight that contains ray $\tau_3$. The ray complementary to $\tau_3$ intersects  $\partial E$ at a point $z_3$. Let us draw a ray $\zeta_3$ starting at $z_3$ that also does not intersect $E$. Among all rays starting at point $x$ and crossing rays $\tau_3$ and $\zeta_3$ there exist two distinct inner supporting rays different from $\tau_3$, which can be shown as in prove of Lemma \ref{lemm4}. This contradicts the fact that $E$ has only two inner supporting rays by our assumption.

 Thus, further we consider sets $E$ that have three and only three inner supporting rays $\xi_i$, $i=1,2,3$, starting at the point $x$. Each point $y_i$ belongs to the boundary of one of two neighboring components. Depending on which boundary of two neighboring components each point $y_i$ belongs, we consider two possible cases: 1) the boundary of each component contains only one point $y_i$  ( $y_i\in\partial E_i$, $i=1,2,3$), Figure 10 a); two points $y_i$ belong to the boundary of the same component ($y_1, y_2\in \partial E_1$, $y_3\in \partial E_3$),  Figure 10 b).
\begin{center}
\includegraphics[width=12 cm]{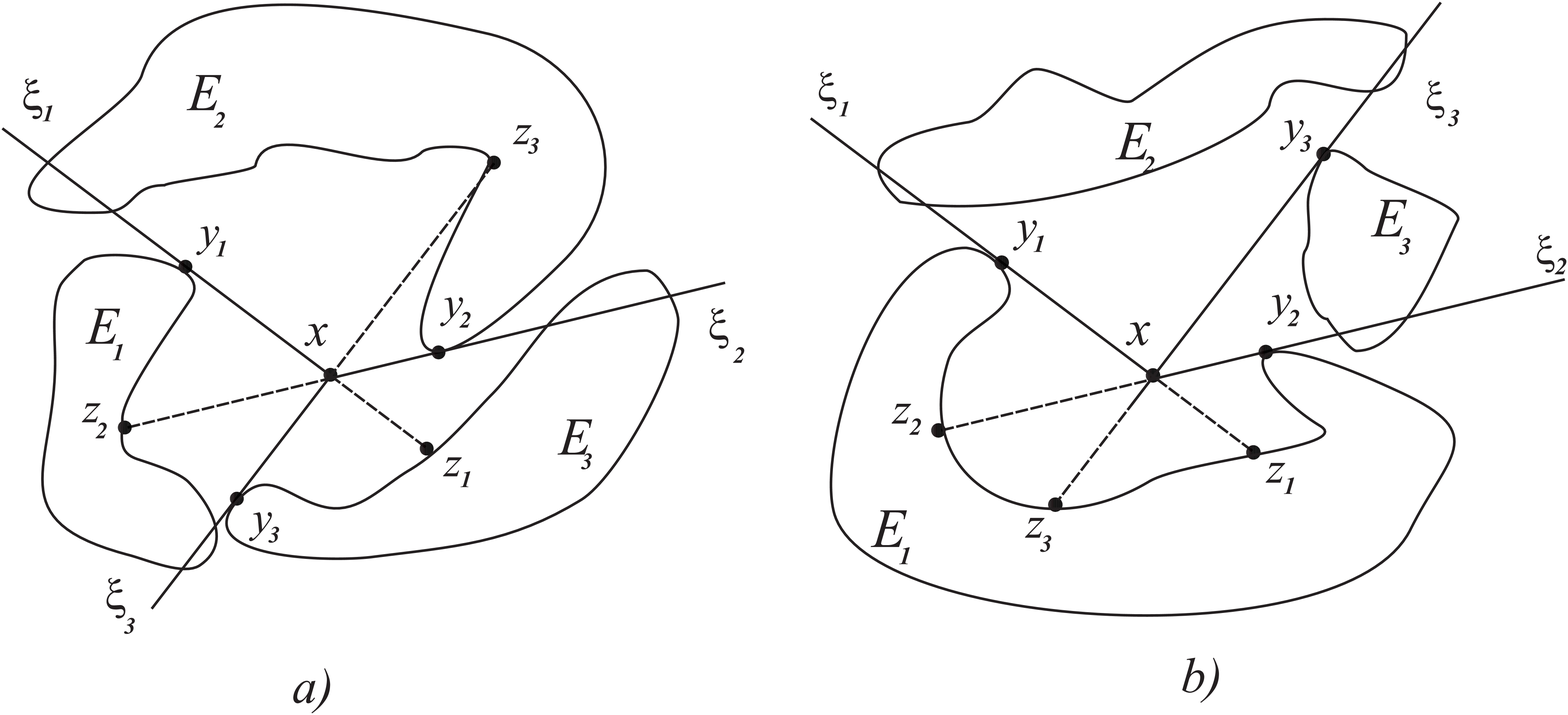}
\end{center}
\small\begin{center} Fig. 10 \end{center}
\normalsize

1) Without loss of generality, let us consider the polar coordinate system where point $x$ is the pole and inner supporting  ray $\xi_1$ is the polar axis. Since $E$ is weakly $1$-semiconvex, there exists a ray $\gamma_1$ starting at the point $y_1$ and not intersecting $E$, Figure 11. The ray $\gamma_1$ should lie above the straight that contains ray $\xi_1$.  Let us chose a positive angular coordinate determined by the angle $\varphi$ between the ray $\xi_1=\xi_1(0)$ and a ray $\xi_1(\varphi)$, $0<\varphi<\pi$, that intersects ray $\gamma_1$.   By the proof of Lemma \ref{lemm4}, there exists an angle $0<\varphi'<\pi$ such that the ray $\xi_1(\varphi')$ is an inner supporting ray and for $\varepsilon>0$, small enough, rays $\xi_1(\varphi)$, $\varphi'<\varphi<\varphi'+\varepsilon$, intersect $E$ before they intersect $\gamma_1$. Suppose  ray $\xi_1(\varphi')$ is tangent to the boundary of one of the components $E_i$, $i=1,2,3$, at a point $y'$.  Let the third inner supporting ray $\xi_1(\varphi'')$ touch $E$ at a point $y''$ and $0<\varphi''<\varphi'$, then points $y_1$, $y'$ belong to the boundary of the component  $E_1$. If $\varphi''>\varphi'$, then one of the points $y_1$, $y'$ belongs to the boundary of the same component that point $y''$ does. Both cases contradict the fact that the boundary of each component contains only one point.
  \begin{center}
\includegraphics[width=8 cm]{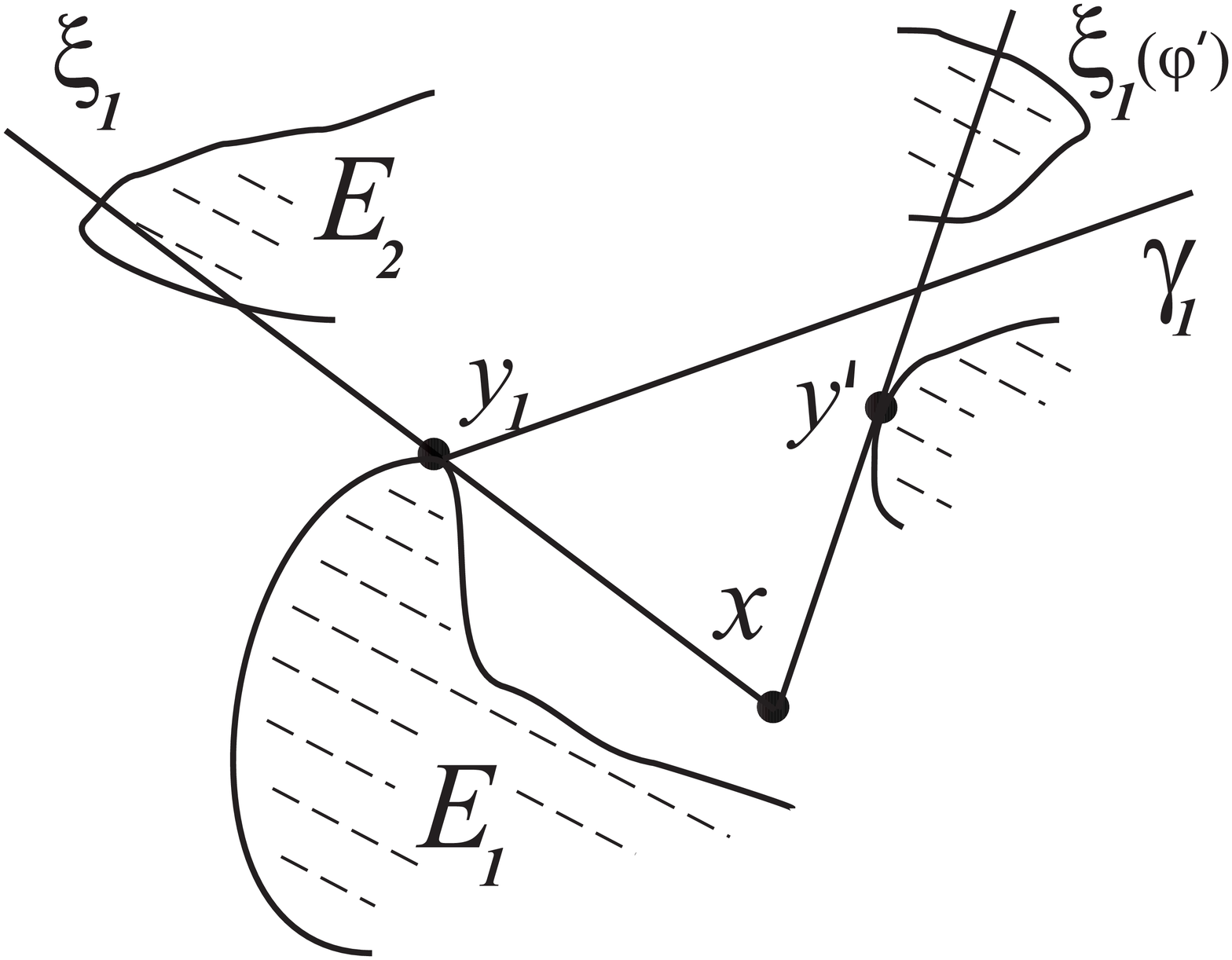}
\end{center}
\small\begin{center} Fig. 11 \end{center}
\normalsize

2)  Let $\alpha$ be the angle of the sector between rays $\xi_1$, $\xi_2$ which contains component $E_1$. Then two cases are possible.

a) $\angle\alpha>\pi$, Figure 12. Then ray $\xi_3$ is contained in the complementary of $\angle\alpha$ to $2\pi$ and point $z_3\in\partial E_1$. Let us draw a ray $\zeta$ that starts at point $z_3$ and does not intersect $E$. We consider the polar coordinate system $(\xi(\varphi),\rho)$  where point $x$ is the pole, the  ray complementary to the inner supporting ray $\xi_3$ is the polar axis $\xi(0)$, and a positive angular coordinate is determined by the angle $\varphi$ between the ray $\xi(0)$ and a ray $\xi(\varphi)$, $0<\varphi<\pi$, that intersects ray $\zeta$. Then there exists an angle $0<\varphi'<\pi$ such that the ray $\xi_4=\xi(\varphi')$ is an inner supporting ray and for $\varepsilon>0$ small enough rays $\xi(\varphi)$, $\varphi'<\varphi<\varphi'+\varepsilon$, intersect $E$ before they intersect ray $\zeta$. Since ray $\zeta$ does not coincide with inner supporting ray  $\xi_3$, ray $\xi_4$ does not coincide with ray  $\xi_3$ too.  Since inner supporting rays $\xi_1$, $\xi_2$ lie in distinct half-planes with respect to the straight that contains ray $\xi_3$, ray $\xi_4$ does not coincide with the ray $\xi_i$, $i\in\{1,2\}$, that lies in the other half-plane. If ray $\xi_4$ coincides with third inner supporting ray, then we have the case when two inner supporting rays coincide with a ray that intersects a unique component, which is not possible too. So, we have four distinct  inner supporting rays, which contradicts the condition that there are three of them.
\begin{center}
\includegraphics[width=8 cm]{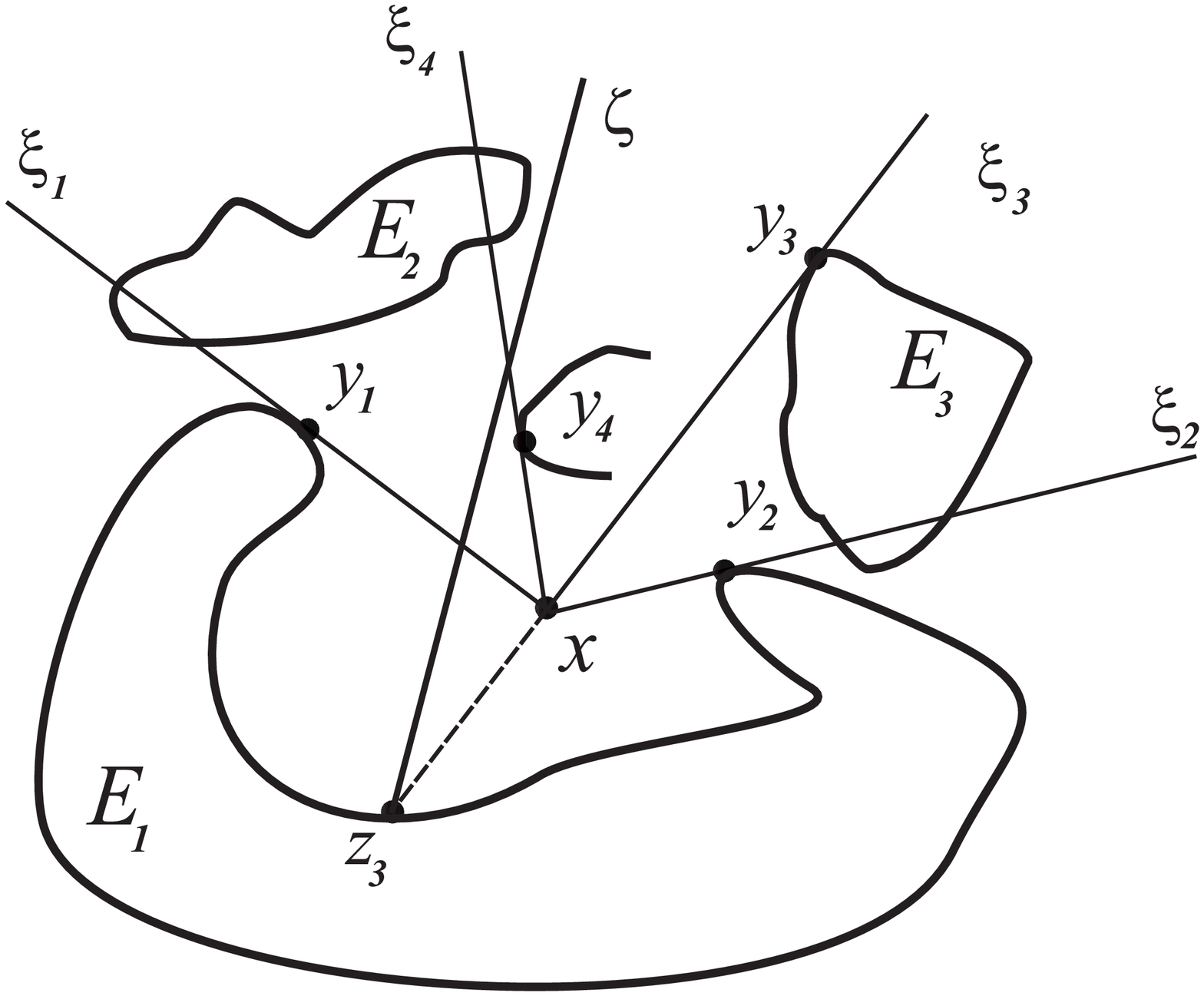}
\end{center}
\small\begin{center} Fig. 12 \end{center}
\normalsize

b) $\angle\alpha\le\pi$.  If the ray complementary to the inner supporting ray $\xi_3$ is contained in the sector $S$ between the rays complementary to $\xi_1$, $\xi_2$, then the statements should be as in case 2 b).

Let $\alpha_1$, $\alpha_2$ be angles between nearest-neighbor rays $\xi_3$, $\xi_1$  and $\xi_3$, $\xi_2$, respectively, Figure 13. Suppose ray $\xi_3$ does not lie in sector $S$, then $\alpha_1\ne\alpha_2$. Let $\alpha_1>\alpha_2$, for definiteness, with that  $\alpha_1>\pi$. A ray that starts at point $y_1$ and does not intersect $E$ should lie above the straight that contains $\xi_1$, by smoothness of $\partial E$ and 1-nonsemiconvexity of $E$ at the point $x$. But any such a ray intersects a part of the component which lies in the sector between rays $\xi_3$, $\xi_1$. Then  $y_1$ is a point of 1-nonsemiconvexity, which contradicts the theorem conditions.
\begin{center}
\includegraphics[width=8 cm]{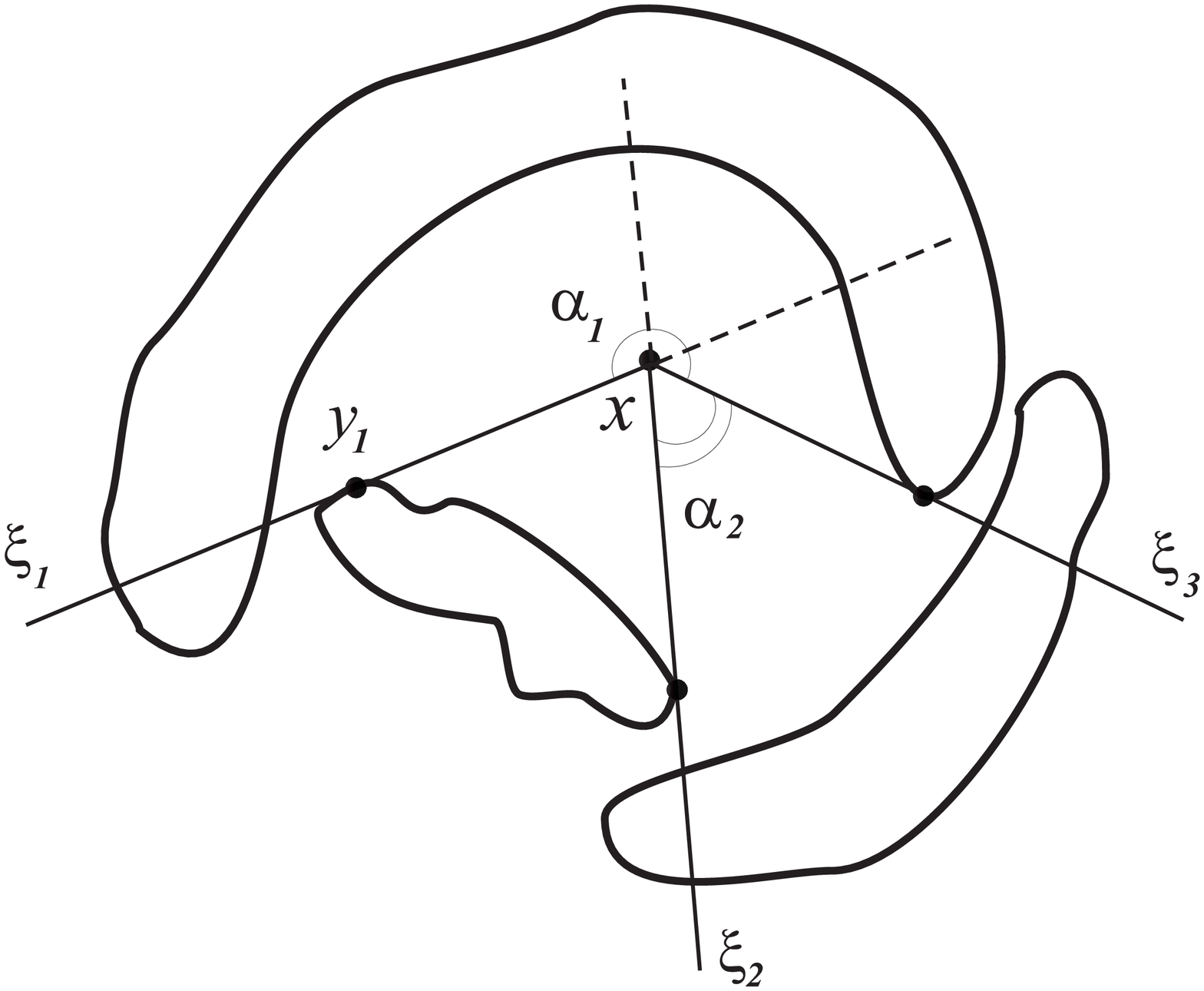}
\end{center}
\small\begin{center} Fig. 13 \end{center}
\normalsize

Thus, set $E$ does not consist of three components. Example \ref{examp3} completes the proof.

\begin{corol} Let an open set $E\subset \mathbb{R}^2$ be weakly $1$-semiconvex but not $1$-semiconvex and consist of four simply connected components with smooth boundary. Then none of its components is projected on the others from a point of $1$-nonsemiconvexity of $E$.
\end{corol}

\noindent{\bf Proof.} Let $E=\bigcup_{i=1}^4E_i$ consist of four simply connected components with smooth boundary and $x\in \mathbb{R}^2\setminus \overline{E}$ be a point of $1$-nonsemiconvexity. Without loss of generality, suppose $E_1$ is projected from $x$  on at least one of the other components. Then the set, consisting only of components $E_i$, $i=2,3,4$, is weakly $1$-semiconvex and not $1$-semiconvex, which contradicts Theorem \ref{theor3}.

\end{document}